\RequirePackage{amsmath}
\documentclass[envcountsect,runningheads]{llncs}
\usepackage{algorithm,algpseudocode}
\usepackage{amsmath,amssymb,amsfonts}

\usepackage{bm}
\usepackage{booktabs}
\usepackage{color}
\usepackage[T1]{fontenc}
\usepackage{graphicx}
\usepackage{hyperref}
\usepackage{mathrsfs,mathtools}
\usepackage{tikz-cd}

\makeatletter
\def\sbar{\accentset{{\cc@style\underline{\mskip13mu}}}}
\def\mbar{\accentset{{\cc@style\underline{\mskip25mu}}}}
\makeatother

\newcommand{\Aut}[1]{\mathrm{Aut}(#1)}

\title{Enumeration algorithm for genus-4 superspecial hyperelliptic curves with automorphism group $\hspace{-0.3mm}\mathbf{Q}_8$}
\titlerunning{Genus-4 superspecial hyperelliptic curves with ${\rm Aut} \cong {\rm Q}_8$}
\author{Takara Taniguchi\inst{1}, Ryo Ohashi\inst{1}, and Tsuyoshi Takagi\inst{1}}
\institute{Graduate School of Information Science and Technology, The University of Tokyo.
\\e-mail: hiroshi-tani@g.ecc.u-tokyo.ac.jp}

\begin{document}
\maketitle
\begin{abstract}
In this paper, we propose an algorithm to enumerate genus-4 superspecial hyperelliptic curves whose automorphism groups isomorphic to the quaternion group.
By implementing this algorithm with Magma, we successfully obtain the number of isomorphism classes of such curves in every characteristic $7 \leq p < 10000$.
Interestingly, the experimental results lead us to the conjecture that there exist exactly $[p/48]$ isomorphism classes of such curves if $p \equiv 1,7 \pmod{8}$, whereas such curves exist if $p \equiv 3,5 \pmod{8}$.
\end{abstract}

\section{Introduction}
Superspecial curves are important research objects in algebraic geometry (for the definition, see Definition \ref{definition}), and their enumeration is a fundamental problem in this field.
The number $N_{g,p}$ of all $\overline{\mathbb{F}}_p$-isomorphism classes of superspecial curves of genus $g$ in characteristic $p>0$ is already known for $g \leq 3$;\,
Deuring~\cite{Deuring1941DieTD} showed that there exists a supersingular elliptic curve (i.e., the case $g=1$) in every characteristic $p \geq 2$, and also determined $N_{1,p}$.
Moreover, the number $N_{2,p}$ (resp. $N_{3,p}$) was determined by Ibukiyama, Katsura, and Oort~\cite{ibukiyama_kasuura_oort} based on the result of Hashimoto and Ibukiyama~\cite{hashimoto1980class} (resp. Hashimoto~\cite{Hashimoto1983ClassNO}).

Unfortunately, the value of $N_{g,p}$ for $g \geq 4$ is unknown in general; \,it is even an open problem whether a genus-$g$ superspecial curve over $\overline{\mathbb{F}}_p$ exists in general $p>0$.
On the other hand, there are several computational results that suggest the existence of a superspecial curve of genus 4;\,
Kudo, Harashita, and Howe~\cite{KHH} showed that there exists a superspecial curve of genus 4 for every $5 \leq p < 20000$ with $p \neq 7$.
Additionally, the existence of a superspecial \emph{hyperelliptic} curve of genus 4 for every $19 \leq p \leq 6691$ is shown by Ohashi, Kudo, and Harashita~\cite{OKH}.
One of their strategies was to consider only on genus-4 curves whose automorphism groups contain the Klein-4 group ${\rm V}_4$.
As another result, Kudo, Nakagawa, and Takagi~\cite{kudo2022efficientsearchsuperspecialhyperelliptic} showed the existence of superspecial genus-4 hyperelliptic curves whose automorphism groups contain the cyclic group ${\rm C}_6$ of order 6 for every $17 \leq p < 1000$ with $p \equiv 2 \pmod{3}$.

This paper serves as a ``sequel'' to the work by Ohashi and Kudo~\cite{OHASHI2024100020klein}, in which they enumerated genus-4 superspecial hyperelliptic curves whose automorphism groups \emph{properly} contain ${\rm V}_4$.
In their paper, the enumeration of genus-4 superspecial hyperelliptic curves with automorphism groups isomorphic to ${\rm Q}_8$ was left as a future work, where ${\rm Q}_8$ denotes the quaternion group.
Hence, we study the number of all $\overline{\mathbb{F}}_p$-isomorphism classes of such curves.
We first show that any genus-4 hyperelliptic curve whose automorphism group contains ${\rm Q}_8$ can be written as
\[
    H_a: y^2 = x(x^4-1)(x^4+ax^2+1)
\]
for an $a \notin \{\pm 2\}$ (see Proposition~\ref{Q8}).
Corollary~\ref{col:num_isomorphism} provides a necessary and sufficient condition on $a$ for $H_a$ to have the automorphism group isomorphic to ${\rm Q}_8$.
Also, we give in Theorem~\ref{thm:autQ8} a necessary and sufficient condition on $a,b$ for $H_a \cong H_b$ holds.
Using these results, we then propose in Section 3.2 an algorithm for enumerating superspecial hyperelliptic curves of genus 4 with automorphism group isomorphic to ${\rm Q}_8$.
As a result of the implementation of our enumeration algorithm in Magma, we found that the number of such curves for every $7 \leq p < 10000$ can be described by the following formula:
\[
    \left\{
    \begin{array}{ll}
        [p/48] & \ \text{ if } p \equiv 1,7 \!\!\pmod{8},\\
        0 & \ \text{ if } p \equiv 3,5 \!\!\pmod{8}.
    \end{array}
    \right.
\]
Specifically, the experimental results for $7 \leq p < 200$ are listed in Table~\ref{tbl:exp_result}.

\section{Preliminaries}
In this section, we recall the definition of hyperelliptic curves, along with the classification of automorphism groups of genus-4 hyperelliptic curves.
Subsequently, we review a criterion for determining whether a hyperelliptic curve is superspecial.
Let $p > 2$ be a prime.
\begin{definition}
A \textbf{hyperelliptic curve} over $\overline{\mathbb{F}}_p$ is a curve of genus $g > 1$ given by the equation
\[
    H: y^2 = f(x), \quad \deg{f(x)} = 2g+1 \text{ or } 2g+2,
\]
where $f(x) \in \overline{\mathbb{F}}_p[x]$ is a square-free polynomial.
\end{definition}
\noindent For a curve $C$ over $\overline{\mathbb{F}}_p$, the set of automorphisms over $\overline{\mathbb{F}}_p$ of $C$ is called the \emph{(full) automorphism group} of $C$, which we denote by ${\rm Aut}(C)$.
The automorphism group of any hyperelliptic curve contains a subgroup of order 2.
\begin{definition}
The \textbf{reduced automorphism group} of a hyperelliptic curve $H$ is defined to be
\[
    {\rm RA}(H) \coloneqq {\rm Aut}(H)/\langle [-1]\rangle
\]
where $[-1]$ denotes the hyperelliptic involution of $H$.
A reduced automorphism of a hyperelliptic curve $H$ can be regarded as an automorphism of $\mathbb{P}^1$.
\end{definition}

In the following, we denote by ${\rm C}_n$ (resp. ${\rm D}_n$) the cyclic (resp. dihedral) group of order $n$ (resp. $2n$) where $n \geq 2$.
We also use ${\rm V}_4$ and ${\rm Q}_8$ to denote the Klein 4-group and the quaternion group, respectively.
The classification of the possible automorphism groups of genus-4 hyperelliptic curves in characteristic $p \geq 7$ is given in the following theorem:

\begin{theorem}[{\cite[Theorem C]{kudo2022arxiv}}]\label{thm:kudo_autq8sadameru}
Assume that $p \geq 7$.
The automorphism group of a genus-4 hyperelliptic curve $H$ over $\overline{\mathbb{F}}_p$ is isomorphic to either of the 12 finite groups in the first column of Table~\ref{tbl:autclassify}.
In each type of ${\rm Aut}(H)$, the curve $H$ is isomorphic to $y^2 = f(x)$ given in the second column of the table where $a,b,c,d \in \overline{\mathbb{F}}_p$.
\end{theorem}

\begin{table}[htbp]
\caption{Classification of hyperelliptic curves of genus 4 over $\overline{\mathbb{F}}_p$ with $p \geq 7$ according to their automorphism groups.}\label{tbl:autclassify}\centering
\begin{tabular}{@{}ll@{}}
\toprule
$\Aut{H}$ & Defining equation $y^2=f(x)$ of $H$  \\ \midrule
        ${\rm C_2}$ & $y^2 = \text{(square-free $f(x)$ of degree $9$ or $10$)}$ \\
        ${\rm V_4}$ & $y^2 = x^{10} + a x^8 + b x^6 + cx^4 + d x^2 +1$ \\ 
        ${\rm C_4}$ & $y^2 = x^{9} + a x^7 + b x^5 + cx^3 + x$ \\ 
        ${\rm C_6}$ & $y^2 = x^{10} + a x^7 + b x^4 + x$ \\[0.5mm]
        ${\rm D_4}$ & $\!\left\{\begin{array}{l} y^2 =  x^9 + a x^7 + b x^5 + a x^3 + x, \,\text{ or }\\ y^2 =  x^{10} + a x^8 + b x^6 + bx^4 + a x^2 + 1\end{array}\right.\vspace{1mm}$\\
        ${\rm Q_8}$ & $y^2 = x(x^4-1)(x^4 + a x^2 + 1) $  \\ 
        ${\rm D_8}$ & $y^2 = x^{9} + a x^5 + x$  \\ 
        ${\rm D_{10}}$ & $y^2 = x^{10} + a x^5 + 1$ \\ 
        $\mathrm{SL}_2(\mathbb{F}_3)$ & $y^2 = x(x^4-1)(x^4 + 2\sqrt{-3} x^2 + 1)$ \\ 
        ${\rm C_{16}} \hspace{-0.1mm}\rtimes\hspace{-0.1mm} {\rm C_2}$ & $y^2 = x^{9} + x$  \\ 
        ${\rm C_{18}}$ & $y^2 = x^{10} + x$ \\ 
        ${\rm C_5} \hspace{-0.1mm}\rtimes\hspace{-0.1mm} {\rm D_4}$ & $y^2 = x^{10} + 1$  \\ 
        \bottomrule
\end{tabular}
\end{table}

Next, we recall the definition of superspecial curves:
\begin{definition}\label{definition}
A \textbf{superspecial curve} is a curve whose Jacobian is isomorphic to a product of supersingular elliptic curves.
\end{definition}


\noindent It is known that a curve is superspecial if and only if its \emph{Cartier-Manin matrix} (see \cite[$\S 2.2$]{AH} for the definition) is the zero matrix.
In particular, the Cartier-Manin matrix of a hyperelliptic curve can be written as in the following theorem:

\begin{theorem}[{\cite[Section 4.1]{YUI1978378}}]
Let $H: y^2=f(x)$ be a genus-$g$ hyperelliptic curve over $\overline{\mathbb{F}}_p$.
The Cartier-Manin matrix of $H$ is given as
    \begin{align*}
        M_H \coloneqq
        \begin{pmatrix} 
            c_{p-1} & c_{p-2} & \dots  & c_{p-g} \\
            c_{2p-1} & c_{2p-2} & \dots  & c_{2p-g} \\
            \vdots & \vdots & \ddots & \vdots \\
            c_{gp-1} & c_{gp-2} & \dots  & c_{gp-g}
        \end{pmatrix}
    \end{align*}
where we write
\[
    f(x)^{(p-1)/2} = \sum_i c_ix^i \text{ with } c_i \in \overline{\mathbb{F}}_p.
\]
\end{theorem}

The superspeciality of certain hyperelliptic curves has already been determined.
We introduce such a example in the following theorem:
\begin{theorem}[{\cite[Theorem 2]{Valentini}}]\label{G32}
Assume that $p \geq 2g+1$.
The genus-$g$ hyperelliptic curve $y^2 = x^{2g+1}-x$ over $\overline{\mathbb{F}}_p$ is superspecial if and only if $p \equiv -1,2g+1 \pmod{4g}$.
\end{theorem}

\section{Proposed method}
In this section, we begin by investigating some properties of genus-4 hyperelliptic curves whose automorphism group contains the quaternion group ${\rm Q_8}$. 
We then propose an algorithm for enumerating the isomorphism classes of superspecial curves among them.

\subsection{A genus-4 hyperelliptic curve whose automorphism group contains ${\rm Q_8}$}

Let $p \geq 7$ be a prime integer, and we consider a genus-4 hyperelliptic curve over $\overline{\mathbb{F}}_p$ defined by the equation
\begin{equation*}
    H_a: y^2 = x(x^4-1)(x^4+ax^2 + 1) \eqqcolon f_a(x)
\end{equation*}
with $a \neq \pm 2$.
Here, the condition $a \neq \pm 2$ is imposed to ensure that $H_a$ is a curve (i.e., non-singular).
This curve has two automorphisms
\begin{align*}
    \sigma: H_a \longrightarrow H_a\,;\ (x,y) &\longmapsto (-x,iy),\\
    \tau: H_a \longrightarrow H_a\,;\ (x,y) &\longmapsto (1/x,iy/x^5)
\end{align*}
with $\sigma^2=\tau^2 = -1$ and $\sigma\tau = -\tau\sigma$.
Hence, the automorphism group of $H_a$ always contains the quaternion group ${\rm Q}_8$.
Conversely, one can show the following:
\begin{proposition}\label{Q8}
Any genus-4 hyperelliptic curve over $\overline{\mathbb{F}}_p$ whose automorphism group contains the quaternion group ${\rm Q}_8$ is isomorphic to $H_a$ for an $a \notin \{2,-2\}$.
\end{proposition}
\begin{proof}
Among the 12 groups listed in Table 1, the ones which contain ${\rm Q}_8$ are only ${\rm Q}_8,\,{\rm SL}_2(\mathbb{F}_3)$, and ${\rm C}_{16} \hspace{-0.4mm}\rtimes {\rm C}_2$.
It immediately follows from Table 1 that a genus-4 hyperelliptic curve $H$ such that ${\rm Aut}(H) \cong {\rm Q}_8$ is isomorphic to $H_a$ for an $a \neq \pm 2$.
Additionally, we have that
\begin{itemize}
\item If ${\rm Aut}(H) \cong {\rm SL}_2(\mathbb{F}_3)$, then $H \cong H_{2\sqrt{3}i}$;
\item If ${\rm Aut}(H) \cong {\rm C}_{16} \hspace{-0.4mm}\rtimes {\rm C}_2$, then $H \cong H_0$;
\end{itemize}
by Table 1 again.
We note that the curve $y^2 = x^9+x$ is isomorphic to $H_0: y^2 = x^9-x$.
\end{proof}\smallskip

Let us describe the explicit structures of the reduced automorphism groups of $H_a$ in each case.
The above automorphisms $\sigma$ and $\tau$ of $H_a$ induces the automorphisms
\[
    \sigma: \mathbb{P}^1 \longrightarrow \mathbb{P}^1\,;\,x \longmapsto -x \ \text{ and }\  \tau: \mathbb{P}^1 \longrightarrow \mathbb{P}^1\,;\,x \longmapsto 1/x
\]
of $\mathbb{P}^1$, respectively (we denote them by $\sigma$ and $\tau$ by abuse of notations).
If ${\rm Aut}(H_a) \cong {\rm Q}_8$, then one can easily check that ${\rm RA}(H_a) = \langle \sigma,\tau \rangle$, which is isomorphic to the Klein-4 group.
Additionally, we have the following:
\begin{itemize}
\item The reduced automorphism group of $H_{2\sqrt{3}i}$ is generated by $\sigma$ and
\[
    \omega: \mathbb{P}^1 \longrightarrow \mathbb{P}^1\,;\ x \longmapsto \frac{x+i}{x-i}.
\]
This is isomorphic to the alternating group of degree $4$.
Note that $\omega$ is of order $3$ and $\omega\sigma\omega^{-1} = \tau$ holds.
A straightforward computation shows that 
\[
    \{\rho(x) \mid \rho \in {\rm RA}(H_{2\sqrt{3}i})\} = \biggl\{\pm x, \pm\frac{1}{x}, \pm\frac{x+i}{x-i}, \pm\frac{x-i}{x+i}, \pm i\frac{x+1}{x-1}, \pm i\frac{x-1}{x+1}\biggr\},
\]
and one can check that the order-2 automorphisms in ${\rm RA}(H_{2\sqrt{3}i})$ are only $\sigma,\tau$, and $\sigma\tau$.
\item The reduced automorphism group of $H_0$ is generated by $\tau$ and
\[
    \mu: \mathbb{P}^1 \longrightarrow \mathbb{P}^1\,;\ x \longmapsto \zeta x,\ \text{ where } \zeta \text{ is a primitive $8$th root of unity }.
\]
This is isomorphic to the dihedral group of order $16$, that is,
\[
    \{\rho(x) \mid \rho \in {\rm RA}(H_0)\} = \{\zeta^kx \mid k = 0,\ldots,7\} \sqcup \{\zeta^k/x \mid k = 0,\ldots,7\}.
\]
One can check that the order-2 automorphisms in ${\rm RA}(H_0)$ are $\mu^4 = \sigma$, or $\mu^k\tau$ with $k \in \{0,\ldots,7\}$.
\end{itemize}

Our aim in this subsection is to provide a necessary and sufficient condition on $a,b$ for $H_b$ to be isomorphic to $H_a$.
For our purposes, we write
\begin{alignat*}{2}
    f_a(x) &= x(x^4-1)(x^2-s^2)(x^2-1/s^2) && \ \text{ with }\ s^2+1/s^2 = -a,\\
    f_b(x) &= x(x^4-1)(x^2-t^2)(x^2-1/t^2) && \ \text{ with }\ t^2+1/t^2 = -b.
\end{alignat*}
and suppose that there exists an isomorphism $\psi: H_b \rightarrow H_a$.
We prepare the following three propositions:
\begin{proposition}
If $\psi: H_b \rightarrow H_a$ satisfies that $\psi(-x) = -\psi(x)$, then $b = \pm a$.
\end{proposition}
\begin{proof}
Write the isomorphism $\psi$ as
\begin{equation*}
    \psi(x) = \frac{\alpha x+\beta}{\gamma x+\delta} \ \text{ with }\, \alpha,\beta,\gamma,\delta \in k.
\end{equation*}
The condition $\psi(-x) = -\psi(x)$ is equivalent to that $2\alpha\gamma x^2 = 2\beta\delta$ for all $x$, which implies that $\alpha\gamma = \beta\delta = 0$.
If $\alpha=\beta=0$ (resp. $\gamma=\delta=0$) holds, then we have $\psi(x) = 0$ (resp. $\psi(x) = \infty$) for all $x$, but this contradicts that $\psi$ is an isomorphism.
Hence, either $\beta=\gamma=0$ or $\alpha=\delta=0$ must hold, and $\psi$ has the form
\[
    \psi(x) = \lambda x\ \text{ or }\,\lambda/x\,\text{ for all }x
\]
where we define $\lambda \coloneqq \alpha/\delta$ (resp. $\lambda \coloneqq \beta/\gamma$) in the former (resp. latter) case.
In each case $\psi(\{0,\infty\}) = \{0,\infty\}$ and since the isomorphism $\psi$ maps the set of branch points of $H_a$ to that of $H_b$, we have that
\[
    \psi(\{\pm 1,\pm i,\pm s, \pm 1/s\}) = \{\pm 1,\pm i,\pm t, \pm 1/t\}.
\]
Additionally, since $\psi(-1) = -\psi(1)$, we have that $\psi(\{\pm 1\}) = \{\pm 1\},\{\pm i\},\{\pm t\}$, or $\{\pm 1/t\}$.
\begin{enumerate}
\item If $\psi(\{\pm 1\}) = \{\pm 1\}$, then we have that $\lambda = \pm 1$.
This leads to that $\psi(\{\pm i\}) = \{\pm i\}$ and
\[
    \{\pm t,\pm 1/t\} = \psi(\{\pm s,\pm 1/s\}) = \{\pm s,\pm 1/s\}.
\]
Since the sum of the squares of the four elements in each set is equal, we have
\begin{align*}
    -2b &= 2(t^2+1/t^2)\\
    &= t^2+(-t)^2+(1/t)^2+(-1/t)^2\\
    &= s^2+(-s)^2+(1/s)^2+(-1/s)^2\\
    &= 2(s^2+1/s^2) = -2a
\end{align*}
which tells us that $b=a$.
\item If $\psi(\{\pm 1\}) = \{\pm i\}$, then we have that $\lambda = \pm i$.
This leads to that $\psi(\{\pm i\}) = \{\pm 1\}$ and
\[
    \{\pm t,\pm 1/t\} = \psi(\{\pm s,\pm 1/s\}) = \{\pm is,\pm i/s\}.
\]
Hence, we have that $2(t^2+1/t^2) = -2(s^2+1/s^2)$ similarly to the case (1), which tells us that $b=-a$.
\item If $\psi(\{\pm 1\}) = \{\pm t\}$, then $\psi(\{\pm i\}) = \{\pm ti\}$ is either of $\{\pm 1\},\{\pm i\}$, or $\{\pm 1/t\}$.
The first (resp. second) case corresponds to $t=\pm i$ (resp. $t=\pm 1$), which contradicts the condition $-(t^2+1/t^2) = b \neq \pm 2$.
We thus have $\{\pm ti\} = \{\pm 1/t\}$, which means that $b = -(t^2+1/t^2) = 0$.
Additionally, it follows from
\[
    \{\pm 1,\pm i\} = \psi(\{\pm s,\pm 1/s\}) = \{\pm t,\pm t/s\}
\]
that $0 = 2t^2(s^2+1/s^2)$ since the sum of the squares of the four elements in each set is equal.
Therefore, we also have that $a = -(s^2+1/s^2) = 0$.
\item If $\psi(\{\pm 1\}) = \{\pm 1/t\}$, then $\psi(\{\pm i\}) = \{\pm i/t\}$ is either of $\{\pm 1\},\{\pm i\}$, or $\{\pm t\}$.
The first (resp. second) case corresponds to $t=\pm i$ (resp. $t=\pm 1$), which contradicts the condition $-(t^2+1/t^2) = b \neq \pm 2$.
We thus have $\{\pm i/t\} = \{\pm t\}$, which means that $b = -(t^2+1/t^2) = 0$.
Additionally, it follows from
\[
    \{\pm 1,\pm i\} = \psi(\{\pm s,\pm 1/s\}) = \{\pm  it,\pm it/s\}
\]
that $0 = -2t^2(s^2+1/s^2)$ similarly to the case (3).
This leads to that $a = -(s^2+1/s^2) = 0$.
\end{enumerate}
In any case (1)--(4), the condition $b = \pm a$ holds, as desired.
\end{proof}

\begin{proposition}
If $\psi: H_b \rightarrow H_a$ satisfies that $\psi(-x) = 1/\psi(x)$, then $b = -2+\frac{16}{a+2}$ or $b = 2+\frac{16}{a-2}$.
\end{proposition}
\begin{proof}
Write the isomorphism $\psi$ as
\begin{equation*}
    \psi(x) = \frac{\alpha x+\beta}{\gamma x+\delta} \ \text{ with }\, \alpha,\beta,\gamma,\delta \in k.
\end{equation*}
The condition $\psi(-x) = 1/\psi(x)$ is equivalent to that $(\alpha^2-\gamma^2)x^2 = \beta^2-\delta^2$ for all $x$, and therefore we obtain that $\gamma=\pm\alpha$ and $\delta=\pm\beta$.
If $\gamma=\alpha$ and $\delta=\beta$ (resp. $\gamma=-\alpha$ and $\delta=-\beta$) holds, then we have that $\psi(x) = 1$ (resp. $\psi(x) = -1$) for all $x$, but this contradicts that $\psi$ is an isomorphism.
Hence $\psi$ has the form
\[
    \psi(x) = \pm\frac{x+\lambda}{x-\lambda}\,\text{ for all }x
\]
where we define $\lambda \coloneqq \alpha/\beta$.
In each case $\psi(\{0,\infty\}) = \{\pm 1\}$ holds, and since the isomorphism $\psi$ maps the set of branch points of $H_a$ to that of $H_b$, we have that
\[
    \psi(\{\pm 1,\pm i,\pm s, \pm 1/s\}) = \{0,\infty,\pm i,\pm t, \pm 1/t\}.
\]
Additionally, since $\psi(-x) = 1/\psi(x)$, we have that $\psi^{-1}(\{0,\infty\}) = \{\pm 1\},\{\pm i\},\{\pm s\}$, or $\{\pm 1/s\}$.
\begin{enumerate}
\item If $\psi^{-1}(\{0,\infty\}) = \{\pm 1\}$, then we have that $\lambda = \pm 1$.
This leads to that $\psi(\{\pm i\}) = \{\pm i\}$ and
\[
    \{\pm t,\pm 1/t\} = \psi(\{\pm s,\pm 1/s\}) = \biggl\{\pm \frac{s+1}{s-1},\pm\frac{s-1}{s+1}\biggr\}.
\]
Since the sum of the squares of the four elements in each set is equal, we have
\begin{align*}
    -2b &= 2(t^2+1/t^2)\\
    &= t^2+(-t)^2+(1/t)^2+(-1/t)^2\\
    &= \biggl(\frac{s+1}{s-1}\biggr)^{\!2}+\biggl(-\frac{s+1}{s-1}\biggr)^{\!2}+\biggl(\frac{s-1}{s+1}\biggr)^{\!2}+\biggl(-\frac{s-1}{s+1}\biggr)^{\!2}\\
    &= 4\frac{s^2+6+1/s^2}{s^2-2+1/s^2} = 4\frac{a-6}{a+2}
\end{align*}
which tells us that $b= -2+\frac{16}{a+2}$.
\item If $\psi^{-1}(\{0,\infty\}) = \{\pm i\}$, then we have that $\lambda = \pm i$.
This leads to that $\psi(\{\pm i\}) = \{\pm 0,\infty\}$ and
\[
    \{\pm t,\pm 1/t\} = \psi(\{\pm s,\pm 1/s\}) = \biggl\{\pm \frac{s+i}{s-i},\pm \frac{s-i}{s+i}\biggr\}.
\]
Hence, we have that $2(t^2+1/t^2) = 4(s^2-6+1/s^2)/(s^2+2+1/s^2)$ similarly to the case (1), which tells us that $b= 2+\frac{16}{a-2}$.
\item If $\psi^{-1}(\{0,\infty\}) = \{\pm s\}$, then $\psi^{-1}(\{\pm i\}) = \{\pm is\}$ must be one of the following: $\hspace{0.3mm}\{\pm 1\},\{\pm i\}$, or $\{\pm 1/s\}$.
The first (resp. second) case implies $s=\pm i$ (resp. $s=\pm 1$), which contradicts $-(s^2+1/s^2) = a \neq \pm 2$.
We thus have $\{\pm is\} = \{\pm 1/s\}$, which means that $a = -(s^2+1/s^2) = 0$.
Additionally, it follows from
\[
    \{\pm t,\pm 1/t\} = \psi(\{\pm 1,\pm i\}) = \biggl\{\pm\frac{1+\zeta}{1-\zeta},\pm\frac{1-\zeta}{1+\zeta}\biggr\}
\]
that $2(t^2+1/t^2) = \pm 12$ since the sum of the squares of the four elements in each set is equal.
Therefore, we obtain that $b = -(t^2+1/t^2) = \pm 6$.
\item If $\psi^{-1}(\{0,\infty\}) = \{\pm 1/s\}$, then $\psi^{-1}(\{\pm i\}) = \{\pm i/s\}$ must be one of the following: $\{\pm 1\},\{\pm i\}$, or $\{\pm s\}$.
The first (resp. second) case implies $s=\pm i$ (resp. $s=\pm 1$), which contradicts $-(s^2+1/s^2) = a \neq \pm 2$.
We thus have $\{\pm i/s\} = \{\pm s\}$, which means that $a = -(s^2+1/s^2) = 0$.
Additionally, it follows from
\[
    \{\pm t,\pm 1/t\} = \psi(\{\pm 1,\pm i\}) = \biggl\{\pm\frac{1+\zeta^{-1}}{1-\zeta^{-1}},\pm\frac{1-\zeta^{-1}}{1+\zeta^{-1}}\biggr\}
\]
that $2(t^2+1/t^2) = \pm 12$ similarly to the case (3).
This leads to that $b = -(t^2+1/t^2) = \pm 6$.
\end{enumerate}
In any case (1)--(4), the condition $b = -2+\frac{16}{a+2}$ or $b = 2+\frac{16}{a-2}$ holds, as desired.
\end{proof}

\begin{proposition}
If $\psi: H_b \rightarrow H_a$ satisfies that $\psi(-x) = -1/\psi(x)$, then $b = 2-\frac{16}{a+2}$ or $b = -2-\frac{16}{a-2}$.
\end{proposition}
\begin{proof}
Write the isomorphism $\psi$ as
\begin{equation*}
    \psi(x) = \frac{\alpha x+\beta}{\gamma x+\delta} \ \text{ with }\, \alpha,\beta,\gamma,\delta \in k.
\end{equation*}
The condition $\psi(-x) = -1/\psi(x)$ is equivalent to that $(\alpha^2-\gamma^2)x^2 = \beta^2-\delta^2$ for all $x$, and therefore we have that $\gamma=\pm i\alpha$ and $\delta=\pm i\beta$.
If $\gamma = i\alpha$ and $\delta=i\beta$ (resp. $\gamma=-i\alpha$ and $\delta=-i\beta$) holds, then we have $\psi(x) = i$ (resp. $\psi(x) = -i$) for all $x$, but this contradicts that $\psi$ is an isomorphism.
Hence $\psi$ has the form
\[
    \psi(x) = \pm i\frac{x+\lambda}{x-\lambda}\,\text{ for all }x
\]
where we define $\lambda \coloneqq \alpha/\beta$.
In each case $\psi(\{0,\infty\}) = \{\pm i\}$ holds, and since the isomorphism $\psi$ maps the set of branch points of $H_a$ to that of $H_b$, we have that
\[
    \psi(\{\pm 1,\pm i,\pm s, \pm 1/s\}) = \{0,\infty,\pm 1,\pm t, \pm 1/t\}.
\]
Additionally, since $\psi(-x) = -1/\psi(x)$, we have that $\psi^{-1}(\{0,\infty\}) = \{\pm 1\},\{\pm i\},\{\pm s\}$, or $\{\pm 1/s\}$.
\begin{enumerate}
\item If $\psi^{-1}(\{0,\infty\}) = \{\pm 1\}$, then we have that $\lambda = \pm 1$.
This leads to that $\psi(\{\pm i\}) = \{\pm 1\}$ and
\[
    \{\pm t,\pm 1/t\} = \psi(\{\pm s,\pm 1/s\}) = \biggl\{\pm i\frac{s+1}{s-1},\pm i\frac{s-1}{s+1}\biggr\}.
\]
Since the sum of the squares of the four elements in each set is equal, we have
\begin{align*}
    -2b &= 2(t^2+1/t^2)\\
    &= t^2+(-t)^2+(1/t)^2+(-1/t)^2\\
    &= \biggl(i\frac{s+1}{s-1}\biggr)^{\!2}\!+\biggl(-i\frac{s+1}{s-1}\biggr)^{\!2}\!+\biggl(i\frac{s-1}{s+1}\biggr)^{\!2}\!+\biggl(-i\frac{s-1}{s+1}\biggr)^{\!2}\\
    &= -4\frac{s^2+6+1/s^2}{s^2-2+1/s^2} = -4\frac{a-6}{a+2}
\end{align*}
which tells us that $b= 2-\frac{16}{a+2}$.
\item If $\psi^{-1}(\{0,\infty\}) = \{\pm i\}$, then we have that $\lambda = \pm i$.
This leads to that $\psi(\{\pm i\}) = \{\pm 0,\infty\}$ and
\[
    \{\pm t,\pm 1/t\} = \psi(\{\pm s,\pm 1/s\}) = \biggl\{\pm i\frac{s+i}{s-i},\pm i\frac{s-i}{s+i}\biggr\}.
\]
Hence, we obtain that $2(t^2+1/t^2) = -4(s^2-6+1/s^2)/(s^2+2+1/s^2)$ similarly to the case (1), which tells us that $b= -2-\frac{16}{a-2}$.
\item If $\psi^{-1}(\{0,\infty\}) = \{\pm s\}$, then $\psi^{-1}(\{\pm 1\}) = \{\pm is\}$ must be one of the following: $\hspace{0.3mm}\{\pm 1\},\{\pm i\}$, or $\{\pm 1/s\}$.
The first (resp. second) case implies $s=\pm i$ (resp. $s=\pm 1$), which contradicts $-(s^2+1/s^2) = a \neq \pm 2$.
We thus have $\{\pm is\} = \{\pm 1/s\}$, which means that $a = -(s^2+1/s^2) = 0$.
Additionally, it follows from
\[
    \{\pm t,\pm 1/t\} = \psi(\{\pm 1,\pm i\}) = \biggl\{\pm i\frac{1+\zeta}{1-\zeta},\pm i\frac{1-\zeta}{1+\zeta}\biggr\}
\]
that $2(t^2+1/t^2) = \pm 12$ since the sum of the squares of the four elements in each set is equal.
Therefore, we obtain that $b = -(t^2+1/t^2) = \pm 6$.
\item If $\psi^{-1}(\{0,\infty\}) = \{\pm 1/s\}$, then $\psi^{-1}(\{\pm 1\}) = \{\pm i/s\}$ must be one of the following: $\{\pm 1\},\{\pm i\},\{\pm s\}$.
The first (resp. second) case implies $s=\pm i$ (resp. $s=\pm 1$), which contradicts $-(s^2+1/s^2) = a \neq \pm 2$.
We thus have $\{\pm i/s\} = \{\pm s\}$, which means that $a = -(s^2+1/s^2) = 0$.
Additionally, it follows from
\[
    \{\pm t,\pm 1/t\} = \psi(\{\pm 1,\pm i\}) = \biggl\{\pm i\frac{1+\zeta^{-1}}{1-\zeta^{-1}},\pm i\frac{1-\zeta^{-1}}{1+\zeta^{-1}}\biggr\}
\]
that $2(t^2+1/t^2) = \pm 12$ similarly to the case (3).
This leads to that $b = -(t^2+1/t^2) = \pm 6$.
\end{enumerate}
In any case (1)--(4), the condition $b = 2-\frac{16}{a+2}$ or $b = -2-\frac{16}{a-2}$ hold, as desired.
\end{proof}

Next, let $\sigma'$ be the order-2 automorphism in ${\rm RA}(H_b)$ such that $\sigma'(x) = -x$ and we consider the following commutative diagram:
\[
    \begin{tikzcd}
        H_a \arrow[d,dashed] \arrow[r,"\psi^{-1}"]& H_b \arrow[d,"\sigma'"]\\
        H_a & H_b \arrow[l,"\psi"]
    \end{tikzcd}
\]
Then, the composition $\psi\sigma'\psi^{-1}$ is an order-2 automorphism in ${\rm RA}(H_a)$.
\begin{lemma}
If the automorphism group of $H_b$ is isomorphic to ${\rm SL}_2(\mathbb{F}_3)$, then $b=\pm 2\sqrt{3}i$.
\end{lemma}
\begin{proof}
Assume that ${\rm Aut}(H_b) \cong {\rm SL}_2(\mathbb{F}_3)$.
It follows from Table 1 that an isomorphism $\psi: H_b \rightarrow H_{2\sqrt{3}i}$ exists.
Since the order-2 automorphisms in ${\rm RA}(H_{2\sqrt{3}i})$ are only $\sigma,\tau$, and $\sigma\tau$, either of the following holds:
\begin{enumerate}
\item[(a)] $\psi\sigma'\psi^{-1} = \sigma$, that is $\psi(-x) = -\psi(x)$ for all $x$.
\item[(b)] $\psi\sigma'\psi^{-1} = \tau$, that is $\psi(-x) = 1/\psi(x)$ for all $x$.
\item[(c)] $\psi\sigma'\psi^{-1} = \sigma\tau$, that is $\psi(-x) = -1/\psi(x)$ for all $x$.
\end{enumerate}
Hence, the above three propositions tell us that $b=\pm 2\sqrt{3}i$.
\end{proof}
\begin{lemma}
If the automorphism group of $H_b$ is isomorphic to ${\rm C}_{16} \hspace{-0.2mm}\rtimes {\rm C}_2$, then $b=0,\pm 6$.
\end{lemma}
\begin{proof}
Assume that ${\rm Aut}(H_b) \cong {\rm C}_{16} \hspace{-0.2mm}\rtimes {\rm C}_2$.
It follows from Table 1 that an isomorphism $\psi: H_b \rightarrow H_0$ exists.
Recall from the structure of ${\rm RA}(H_a)$ that either of the following holds:
\begin{enumerate}
\item[(a)] $\psi\sigma'\psi^{-1} = \sigma$, that is $\psi(-x) = -\psi(x)$ for all $x$.
\item[(b)] $\psi\sigma'\psi^{-1} = \zeta^k\tau$ with $k \in \{0,\ldots,7\}$, that is $\psi(-x) = \zeta^k/\psi(x)$ for all $x$.
\end{enumerate}
In the case (a), we have that $b = 0$ thanks to Proposition 3.2.
The case (b) with $k \in \{1,3,5,7\}$ cannot occur.
Indeed, substituting $x=0$ yields $\psi(0)^2 = \zeta^k$, which contradicts that $\psi$ maps $0$ to a branch point of $H_0$.
Also, if $\psi$ satisfies $\psi\sigma'\psi^{-1} = \zeta^k\tau$ with $k \in \{0,2,4,6\}$, then $\tilde\psi \coloneqq \zeta^{-k/2}\psi$ is an isomorphism satisfying $\tilde\psi\sigma'\tilde\psi^{-1} = \tau$.
Hence we have $b=\pm 6$ thanks to Proposition 3.3.
\end{proof}

\begin{theorem}\label{thm:autQ8}
The curve $H_b$ is isomorphic to $H_a$ if and only if
\begin{equation}
    b \in \biggl\{\pm a,\ \pm\biggl(2-\frac{16}{a+2}\biggr),\ \pm\biggl(2+\frac{16}{a-2}\biggr)\biggr\}.
\end{equation}
\end{theorem}
\begin{proof}
We first show the ``if'' part:
\begin{itemize}
\item If $b=a$, then we can construct an isomorphism $\psi: H_b \rightarrow H_a$ with $\psi(x) = x$.
\item If $b=-a$, then we can construct an isomorphism $\psi: H_b \rightarrow H_a$ with $\psi(x) = ix$.
\item If $b=2-\frac{16}{a+2}$, then we can construct an isomorphism $\psi: H_b \rightarrow H_a$ with $\psi(x) = i\frac{x+1}{x-1}$.
\item If $b=-2+\frac{16}{a+2}$, then we can construct an isomorphism $\psi: H_b \rightarrow H_a$ with $\psi(x) = \frac{x+1}{x-1}$.
\item If $b=2+\frac{16}{a-2}$, then we can construct an isomorphism $\psi: H_b \rightarrow H_a$ with $\psi(x) = \frac{x+i}{x-i}$.
\item If $b=-2-\frac{16}{a-2}$, then we can construct an isomorphism $\psi: H_b \rightarrow H_a$ with $\psi(x) = i\frac{x+i}{x-i}$.
\end{itemize}
Next, we show the ``only if'' part.
\begin{itemize}
\item If ${\rm Aut}(H_a) \cong {\rm SL}_2(\mathbb{F}_3)$, then Lemma 3.5 implies that $a,b \in\hspace{-0.2mm} \{\pm2\sqrt{3}i\}$, which satisfies the equation (3.1).
\item If ${\rm Aut}(H_a) \cong {\rm C}_{16} \hspace{-0.2mm}\rtimes {\rm C}_2$, then Lemma 3.6 implies that $a,b \in \{0,\pm 6\}$, which satisfies the equation (3.1).
\item If ${\rm Aut}(H_a) \cong {\rm Q}_8$, then any isomorphism $\psi: H_b \rightarrow H_a$ satisfies $\psi\sigma'\psi^{-1} = \sigma,\tau$, or $\sigma\tau$. This implies that the equality $\psi(-x) = -\psi(x),1/\psi(x)$, or $-1/\psi(x)$ holds for all $x$.
Therefore, Propositions 3.2--3.4 tells us that $b$ satisfies the equality (3.1).
\end{itemize}
From the above discussions, we obtain this theorem.
\end{proof}

\begin{corollary}\label{col:num_isomorphism}
The automorphism group of $H_a$ with $a \notin \{\pm2\}$ can be determined as follows:
\begin{enumerate}
\item If $a^2 = -12$, then we have $\Aut{H_a} \cong {\rm SL}_2(\mathbb{F}_3)$.
\item If $a \in \{0,\pm 6\}$, then we have $\Aut{H_a} \cong {\rm C}_{16} \hspace{-0.4mm}\rtimes {\rm C}_2$.
\item Otherwise, we have $\Aut{H_a} \cong {\rm Q}_8$.
\end{enumerate}
\end{corollary}\medskip

Next, define $e \coloneqq (p-1)/2$ and
\[
    f_a(x)^e = \sum_{k=0}^{9e}\gamma_k x^k \text{ with } \gamma_k \in \mathbb{F}_p[a].
\]
Let us describe the Cartier-Manin matrix of $H_a$.
\begin{proposition}\label{sign-equality}
For each $k \in \{e,\ldots,9e\}$, we have
\[
    \gamma_{5p-5-k} = (-1)^e\gamma_k.
\]
\end{proposition}
\begin{proof}
Write
\[
    (x^4-1)^e = \sum_{i=0}^{4e}\alpha_ix^i,\ \ \{x(x^4+ax^2+1)\}^e = \sum_{j=e}^{5e}\beta_jx^j.
\]
It follows the binomial theorem that
\[
    \alpha_i = \left\{
    \begin{array}{ll}
    (-1)^{e-i/4}\binom{e}{i/4} & \text{ if $i$ is a multiple of $4$},\\
    0 & \text{ otherwise}
    \end{array}
    \right.
\]
and
\[
    \alpha_{4e-i} = \left\{
    \begin{array}{ll}
    (-1)^{i/4}\binom{e}{e-i/4} & \text{ if $i$ is a multiple of $4$},\\
    0 & \text{ otherwise}.
    \end{array}
    \right.
\]
This tells us that $\alpha_{4e-i} = (-1)^e\alpha_i$ for all $i \in \{0,\ldots,4e\}$.
Moreover, since $x(x^4+ax^2+1)$ is a symmetric polynomial, we see that $\beta_{6e-j} = \beta_j$ for all $j \in \{e,\ldots,5e\}$.
Here, since $f_a(x)^e = (x^4-1)^e\{x(x^4+ax^2+1)\}^e$ is equal to
\[
    \Biggl(\raisebox{0.5mm}{$\displaystyle\sum_{i=0}^{4e}\alpha_ix^i$}\Biggr)\!\Biggl(\raisebox{0.5mm}{$\displaystyle\sum_{j=e}^{5e}\beta_jx^j$}\Biggr) = \sum_{k=e}^{9e}\Biggl(\raisebox{0.5mm}{$\displaystyle\sum_{i+j=k}\alpha_i\beta_j$}\Biggr)x^k,
\]
we obtain that
\[
    \gamma_k = \sum_{i+j=k}\alpha_i\beta_j,\ \ \gamma_{10e-k} = \sum_{i+j=k}\alpha_{4e-i}\beta_{6e-j}
\]
for all $k \in \{e,\ldots,9e\}$.
We finally have the equality
\[
    \gamma_{10e-k} = (-1)^e\sum_{i+j=k}\alpha_i\beta_j = (-1)^e\gamma_k,
\]
as desired.
\end{proof}
\begin{corollary}\label{prop:cmmatrix_fourshape}
For each $a \notin \{\pm2\}$, let $M_a$ be the Cartier-Manin matrix of $H_a$.
Then, we have the following:
\begin{enumerate}
\item If $p \equiv 1 \pmod{4}$, then
\[
    M_a = \begin{pmatrix}
            \gamma_{p-1} & \!0 & \!\gamma_{p-3} & \!0 \\
            0 & \!\gamma_{2p-2} & \!0 & \!\gamma_{2p-4} \\
            \gamma_{2p-4} & \!0 & \!\gamma_{2p-2} & \!0 \\
            0 & \!\gamma_{p-3} & \!0 & \!\gamma_{p-1}
        \end{pmatrix}.
\]
\item If $p \equiv 3 \pmod{4}$, then
\[
    M_a = \begin{pmatrix}
            0 & \!\!\!\gamma_{p-2} & \!\!\!0 & \!\!\!\gamma_{p-4} \\
            \gamma_{2p-1} & \!\!\!0 & \!\!\!\gamma_{2p-3} & \!\!\!0 \\
            0 & \!\!\!-\gamma_{2p-3} & \!\!\!0 & \!\!\!-\gamma_{2p-1} \\
            -\gamma_{p-4} & \!\!\!0 & \!\!\!-\gamma_{p-2} & \!\!\!0
        \end{pmatrix}.
\]
\end{enumerate}
\end{corollary}
\begin{proof}
The polynomial $(x^4-1)^e(x^4+ax^2+1)^e$ has only terms of even degree, and thus we obtain that
\begin{itemize}
\item[\emph{(1)}] If $e$ is even, then $\gamma_k = 0$ for all odd $k$;
\item[\emph{(2)}] If $e$ is odd, then $\gamma_k = 0$ for all even $k$.
\end{itemize}
Then, this corollary follows from Proposition \ref{sign-equality}.
\end{proof}

Define the polynomial ${\sf gcdall} \in \mathbb{F}_p[a]$ as
\[
    \left\{
        \begin{array}{ll}
            {\rm gcd}(\gamma_{p-1},\gamma_{p-3},\gamma_{2p-2},\gamma_{2p-4}) & \text{ if } p \equiv 1 \!\!\pmod{4},\\
            {\rm gcd}(\gamma_{p-2},\gamma_{p-4},\gamma_{2p-1},\gamma_{2p-3}) & \text{ if } p \equiv 3 \!\!\pmod{4}.
        \end{array}
    \right.
\]
Then, it follows from Corollary \ref{prop:cmmatrix_fourshape} that $H_{a_0}$ with $a_0 \neq \pm 2$ is superspecial if and only if ${\sf gcdall}(a_0) = 0$.

\subsection{Enumeration algorithm}\label{algorithm}
First, we explain a method to compute $\textsf{gcdall}$ defined in the end of Section 3.1 through a simple improvement.
Define $F_a(X) \coloneqq (X^2-1)(X^2+aX+1)$ and 
\[
    F_a(X)^e = \sum_{k=0}^{4e}\delta_k X^k \text{ with } \delta_k \in \mathbb{F}_p[a].
\]
A straightforward computation leads to that $\delta_k = \gamma_{2k+e}$ for each $k \in \{0,\ldots,4e\}$.
Using this fact, we can compute $\textsf{gcdall} \in \mathbb{F}_p[a]$ as in Algorithm \ref{alg:calculate_cm}.

\begin{algorithm}[htbp]
    \caption{}\label{alg:calculate_cm}
     \renewcommand{\algorithmicrequire}{\textbf{Input:}}
     \renewcommand{\algorithmicensure}{\textbf{Output:}}
    \begin{algorithmic}[1]
        \Require a prime $p \ge 7$
        \Ensure the polynomial $\textsf{gcdall} \in \mathbb{F}_p[a]$ defined as above
        \State \label{line:calc_cm_end}$F_a(X) \leftarrow (X^2-1)(X^2+aX+1)$
        \If{$p \equiv 1 \pmod{4}$}
            \State \textbf{return} ${\rm gcd}(\delta_{(p-1)/4},\delta_{(p-5)/4},\delta_{(3p-3)/4},\delta_{(3p-7)/4})$
        \Else
            \State \textbf{return} ${\rm gcd}(\delta_{(p-3)/4},\delta_{(p-7)/4},\delta_{(3p-1)/4},\delta_{(3p-5)/4})$
        \EndIf
    \end{algorithmic}
\end{algorithm}


Next, we explain a method to enumerate isomorphism classes of superspecial curves $H_a$.
We suppose that
\begin{equation}\label{assertion}\tag{$\star$}
    \left\{
    \begin{array}{ll}
        \textsf{gcdall}(2) \neq 0,\ \textsf{gcdall}(-2) \neq 0, \text{ and}\\
        {\rm gcd}\bigl(\textsf{gcdall}(a),\textstyle\frac{d}{da}\textsf{gcdall}(a)\bigr) = 1.
    \end{array}
    \right.
\end{equation}
Under these conditions, the polynomial $\mathsf{gcdall}(a)$ has no root at $a = \pm 2$, and all of its roots are simple.
Thus, the degree of $\mathsf{gcdall}$ is equal to the number of $a$ such that $H_a$ is superspecial.
Among them, the cases corresponding to $\operatorname{Aut}(H_a) \supsetneq Q_8$ can be detected as follows:



\begin{itemize}
\item If $\textsf{gcdall}(2\sqrt{-3}) = 0$, then the two curves in case (1) of Corollary \ref{col:num_isomorphism} are superspecial;
\item If $\textsf{gcdall}(0) = 0$, then the three curves in case (2) of Corollary \ref{col:num_isomorphism} are superspecial;
\end{itemize}

\noindent All remaining values of $a$ correspond to ${\rm Aut}(H_a) \cong {\rm Q}_8$, and each such $a$ yields exactly $6$ isomorphic curves from Theorem \ref{thm:autQ8}.
Hence, we have the following algorithm:

\begin{algorithm}[htbp]
    \caption{}
\label{alg:count_automorphism}
     \renewcommand{\algorithmicrequire}{\textbf{Input:}}
     \renewcommand{\algorithmicensure}{\textbf{Output:}}
    \begin{algorithmic}[1]   
        \Require a prime $p \ge 7$
        \Ensure the number of isomorphism classes of genus-4 superspecial hyperelliptic curves with automorphism group $\hspace{-0.3mm}{\rm Q}_8$ in characteristic $p$
        \State Compute $\textsf{gcdall} \in \mathbb{F}_p[a]$ using \textbf{Algorithm 1}
        \State \textbf{if} (\ref{assertion}) does not hold \textbf{then} \,\Return $\perp$ \,\textbf{end if}
        \State $d \leftarrow \deg{(\textsf{gcdall})}$
        \State \textbf{if} $\textsf{gcdall}(2\sqrt{-3}) = 0$ \textbf{then} \,$d \leftarrow d-2$ \,\textbf{end if}
        \State \textbf{if} $\textsf{gcdall}(0) = 0$ \textbf{then}\, $d \leftarrow d-3$ \,\textbf{end if}
        \State \Return $d/6$  
\end{algorithmic}
\end{algorithm}

\begin{remark}
The condition (\ref{assertion}) has been checked to hold for every prime $p$ with $7 \leq p < 10000$.
\end{remark}

\section{Experimental results}
We have implemented and executed our proposed algorithms in Section \ref{algorithm} on Magma\,V2.28-4, using a computer with an AMD EPYC 7742 CPU and 2TB of RAM.
Our code is available at\smallskip
\begin{center}
    \url{https://github.com/taniguchitakara/EnmAlgG4Q8}.\medskip
\end{center}

\noindent As a result, we succeeded in determining the number of isomorphism classes of genus-4 superspecial hyperelliptic curves with automorphism group $\supset {\rm Q}_8$ for all primes $p$ with $7 \leq p < 10000$.
It took a total of 47,856 seconds to obtain these results.
In particular, the results for $p < 200$ are presented in Table~\ref{tbl:exp_result} (we respectively denote ${\rm SL}_2(\mathbb{F}_3)$ and ${\rm C}_{16} \rtimes {\rm C}_2$ by ${\rm G}_{24}$ and ${\rm G}_{32}$).

\begin{table}[htbp]
        \caption{Enumeration result of genus-4 superspecial hyperelliptic curves $H$ such that ${\rm Aut}(H) \supset {\rm Q}_8$, and their classification by types of automorphism groups}\label{tbl:exp_result}\smallskip
        \centering
        \begin{tabular}{@{}ccccc@{}}
        \toprule
        $p$ & ${\rm Q_8}$ & ${\rm G}_{24}$ & ${\rm G}_{32}$ \\  \midrule
         7 & 0 & 0 & 0 \\ 
        11 & 0 & 0 & 0 \\ 
        13 & 0 & 0 & 0 \\ 
        17 & 0 & 0 & 0 \\ 
        19 & 0 & 0 & 0 \\ 
        23 & 0 & 0 & 1 \\ 
        29 & 0 & 0 & 0 \\ 
        31 & 0 & 1 & 0 \\ 
        37 & 0 & 0 & 0 \\ 
        41 & 0 & 1 & 1 \\ 
        43 & 0 & 0 & 0 \\ 
        47 & 0 & 1 & 1 \\ 
        53 & 0 & 0 & 0 \\ 
        59 & 0 & 0 & 0 \\ 
        61 & 0 & 0 & 0 \\ 
        67 & 0 & 0 & 0 \\ 
        71 & 1 & 0 & 1 \\ 
        73 & 1 & 1 & 0 \\ 
        79 & 1 & 1 & 0 \\ 
        83 & 0 & 0 & 0 \\ 
        89 & 1 & 1 & 1 \\ 
        97 & 2 & 0 & 0 \\   
        \bottomrule
        \end{tabular}\hspace{3mm}
        \begin{tabular}{@{}ccccc@{}}
        \toprule
        $p$ & ${\rm Q_8}$ & ${\rm G}_{24}$ & ${\rm G}_{32}$ \\  \midrule
        101 & 0 & 0 & 0 \\ 
        103 & 2 & 0 & 0 \\ 
        107 & 0 & 0 & 0 \\ 
        109 & 0 & 0 & 0 \\ 
        113 & 2 & 0 & 1 \\ 
        127 & 2 & 1 & 0 \\ 
        131 & 0 & 0 & 0 \\ 
        137 & 2 & 1 & 1 \\ 
        139 & 0 & 0 & 0 \\ 
        149 & 0 & 0 & 0 \\ 
        151 & 3 & 0 & 0 \\ 
        157 & 0 & 0 & 0 \\ 
        163 & 0 & 0 & 0 \\ 
        167 & 3 & 0 & 1 \\ 
        173 & 0 & 0 & 0 \\ 
        179 & 0 & 0 & 0 \\ 
        181 & 0 & 0 & 0 \\ 
        191 & 3 & 1 & 1 \\ 
        193 & 4 & 0 & 0 \\ 
        197 & 0 & 0 & 0 \\ 
        199 & 4 & 0 & 0 \\\\
        \bottomrule
        \end{tabular}
\end{table}

From the above experimental result, we finally obtain the following theorem.
\begin{theorem}\label{main}
For each prime $7 \leq p < 10000$, the following assertions are true:
\begin{enumerate}
\item The number of isomorphism classes of superspecial hyperelliptic curves of genus 4 whose automorphism group $\cong {\rm Q}_8$ is equal to
\[
    \left\{
    \begin{array}{ll}
        [p/48] &\ \text{ if } p \equiv 1,7 \!\!\pmod{8},\\
        0 &\ \text{ if } p \equiv 3,5 \!\!\pmod{8}.
    \end{array}
    \right.
\]
Here, brackets $[\,\cdot\,]$ denote the floor function.
\item A hyperelliptic curve of genus 4 with automorphism group $\cong {\rm SL}_2(\mathbb{F}_3)$ is superspecial if and only if $p \equiv 17 \text{ or } 23 \pmod{24}$.
\item A hyperelliptic curve of genus 4 with automorphism group $\cong {\rm C}_{16} \hspace{-0.2mm}\rtimes\hspace{-0.2mm} {\rm C}_2$ is superspecial if and only if $p \equiv 9 \text{ or } 15 \pmod{16}$.
\end{enumerate}
\end{theorem}

By recalling that a hyperelliptic curve of genus 4 with automorphism group $\cong {\rm C}_{16} \hspace{-0.2mm}\rtimes\hspace{-0.2mm} {\rm C}_2$ is isomorphic to
\[
    H_0: y^2 = x^9-x,
\]
we conclude that the statement (3) in Theorem \ref{main} holds for all $p \geq 11$ applying Theorem \ref{G32} in the case $g=4$.
On the other hand, we conjecture that statements (1)(2) in Theorem \ref{main} also hold for general $p \geq 7$.
Theoretical proofs of these conjectures remain as future work.

\section*{Acknowledgments}
This work was supported by JSPS Grant-in-Aid for Young Scientists 25K17225.
This work was also supported by JST CREST Grant Number JPMJCR2113, Japan.

\end{document}